# SUBCRITICAL REGIMES IN THE POISSON BOOLEAN MODEL OF CONTINUUM PERCOLATION

By Jean-Baptiste Gouéré

*Université d'Orléans*

We consider the Poisson Boolean model of continuum percolation. We show that there is a subcritical phase if and only if $E(R^d)$ is finite, where $R$ denotes the radius of the balls around Poisson points and $d$ denotes the dimension. We also give related results concerning the integrability of the diameter of subcritical clusters.

**1. Introduction.** We consider the Poisson Boolean model of continuum percolation. At each point of a homogeneous Poisson point process on the Euclidean space $\mathbb{R}^d$, we center a ball of random radius. We assume that the radii of the balls are independent, identically distributed and independent of the point process. We denote by $\Sigma$ the union of the balls and by $S$ the connected component of $\Sigma$ that contains the origin. In this paper, we are interested in certain properties of $S$ when the density $\lambda$ of the Poisson point process is small.

Let $R$ be one of the random radii. In [2] (see also [3] and [5]), Hall proved that if $E(R^{2d-1})$ is finite, then the set $S$ is almost surely bounded for small enough $\lambda$. If $E(R^d)$ is infinite, then such behavior does not occur: whatever the value of the density $\lambda$, the set $\Sigma$ is almost surely the whole space. In fact, a much more general statement is known. In the book by Meester and Roy ([5], Proposition 7.3), it is proved that in any almost surely nonempty and stationary point process, it is the case that if $E(R^d)$ is infinite, then the whole space is covered almost surely. In this paper, we prove that the set $S$ is almost surely bounded for small enough $\lambda$ if and only if $E(R^d)$ is finite.

Let us denote by $N$ the number of balls contained in $S$. In [2] (see also [3] and [5]), Hall also proved that $N$ is integrable for small enough $\lambda$ if and only if $E(R^{2d})$ is finite. More generally, in [1], Błaszczyszyn, Rau and Schmidt proved, among other things, that for all integer $k \geq 1$, $N^k$ is integrable for small enough $\lambda$ if and only if the moment $E(R^{d(1+k)})$ is finite. In this 4paper, we prove a related result. Let us denote by $D$ the Euclidean diameter of $S$. We prove that for all $s > 0$, $D^s$ is integrable for small enough $\lambda$ if and only if the moment $E(R^{d+s})$ is finite. This is an improvement of a result of Menshikov, Popov and Vachkovskaia [6] which can be stated as follows: if $E(R^{2sd})$ is finite for a positive integer $s$, then $D^s$ is integrable for small enough $\lambda$. This sufficient condition for the finiteness of moments of $D$ was used in [6] to give sufficient conditions for nonpercolation results in a multiscale percolation model.

In all of these results, the main difficulty lies in proving that $S$ is small for small enough $\lambda$ when $R$ is sufficiently integrable. (Proving that this behavior does not occur when $R$ is not sufficiently integrable is easy.) The proofs of [1], [2] and [6] all rely on the introduction of a multitype branching process that dominates the percolation process. Roughly speaking, the branching process is defined so that the size of the first generation dominates the number of balls that contain the origin, the size of the second generation dominates the number of balls that intersect the previous balls, and so on. In this paper, the proof relies on an estimate that can be roughly described as follows (see Proposition 3.1 or Lemma 3.3 for a more precise statement): the probability of $D$ being greater than a real $\alpha$ is bounded above by the square of the probability of $D$ being greater than $\alpha/10$, up to error terms that are due to the existence of large balls.

**2. Notation and statement of the main results.** For the whole of the paper, we fix an integer $d \geq 1$. Let $|\cdot|$ denote the Lebesgue measure on $\mathbb{R}^d$. We denote by $\|\cdot\|$ the Euclidean norm on $\mathbb{R}^d$ and by $B(x,r)$ the open Euclidean ball centered at $x \in \mathbb{R}^d$ with radius $r \geq 0$.

Let $\lambda > 0$ and let $\nu$ be a probability measure on $]0, +\infty[$. Let $\xi$ be a Poisson point process on $\mathbb{R}^d \times ]0, +\infty[$ whose intensity measure is the product of $\lambda|\cdot|$ and $\nu$. We denote by $P_{\lambda,\nu}$, $E_{\lambda,\nu}$ the associated probability measure and expectation, respectively. As distinct points of $\xi$ have distinct coordinates on $\mathbb{R}^d$, we can write

$$\xi = \{(c, r(c)), c \in \chi\},$$

where $\chi$ denotes the projection of $\xi$ on $\mathbb{R}^d$. Note that $\chi$ is a Poisson point process on $\mathbb{R}^d$ whose intensity is $\lambda|\cdot|$. Let us recall that if we condition on $\chi$, then the $r(c), c \in \chi$, are i.i.d. with common distribution $\nu$. We refer to [4, 7, 8] for background on point processes and to [3, 5] for Boolean models.

We are interested in the properties of the following random set:

$$(1) \qquad \Sigma = \bigcup_{c \in \chi} B(c, r(c)).$$

Let $S$ denote the connected component of $\Sigma$ which contains 0. (We let $S = \varnothing$ if 0 does not belong to $\Sigma$.) We say that percolation occurs if the set $S$ is

unbounded:

$$\{\text{percolation}\} = \{S \text{ is unbounded}\}.$$

We prove the following result.

THEOREM 2.1. *For all probability measures $\nu$ on $]0, +\infty[$, the following assertions are equivalent:*

1. *the moment $\int r^d \nu(dr)$ is finite;*
2. *there exists $\lambda_0 > 0$ such that $P_{\lambda,\nu}(Percolation) = 0$ for all positive real $\lambda < \lambda_0$.*

*Moreover, there exists a positive constant $C$ that depends only on the dimension $d$ such that, for all probability measures $\nu$ on $]0, +\infty[$, if the previous assertions hold, then*

$$P_{\lambda,\nu}(Percolation) = 0 \quad \text{for all positive real } \lambda < C\left(\int r^d \nu(dr)\right)^{-1}.$$

Let $D$ denote the Euclidean diameter of $S$:

(2) $$D = \sup_{x,y \in S} \|x - y\|.$$

(We let $D = 0$ if $S$ is empty.) We prove the following result.

THEOREM 2.2. *Let $s$ be a positive real. For all probability measures $\nu$ on $]0, +\infty[$, the following assertions are equivalent:*

1. *the moment $\int r^{d+s} \nu(dr)$ is finite;*
2. *there exists $\lambda_0 > 0$ such that $E_{\lambda,\nu}(D^s)$ is finite for all positive real numbers $\lambda < \lambda_0$.*

*Moreover, there exists a positive constant $C$ that depends only on the dimension $d$ such that, for all positive real $s$ and all probability measures $\nu$ on $]0, +\infty[$, if the previous assertions hold, then*

$$E_{\lambda,\nu}(D^s) \text{ is finite for all positive real numbers } \lambda < C\left(\int r^d \nu(dr)\right)^{-1}.$$

**3. Proofs.** In order to simplify the notation, in the whole of this section, once a probability measure $\nu$ on $]0, +\infty[$ is given, we will denote by $R$ a random variable whose law is $\nu$. Theorems 2.1 and 2.2 are consequences of Lemma 3.8 in Section 3.2 and Lemma 3.9 in Section 3.3.

3.1. *Proof of some inequalities.* In all of this subsection, we fix $\lambda > 0$ and a probability measure $\nu$ on $]0, +\infty[$. In order to simplify the notation, we drop the subscript $\{\lambda, \nu\}$ from the probability measure $P$ and from the expectation symbol $E$. For all Borel subsets $A \subset \mathbb{R}^d$, we define a set $\Sigma(A)$ as follows:

$$\Sigma(A) = \bigcup_{c \in \chi \cap A} B(c, r(c)). \tag{3}$$

We will study percolation through a family of events defined as follows. If $\alpha > 0$ is a real and $x$ is a point of $\mathbb{R}^d$, we say that $G(x, \alpha)$ occurs if the connected component of

$$\Sigma(B(x, 10\alpha)) \cup B(x, \alpha)$$

containing $x$ is not contained in $B(x, 8\alpha)$. By stationarity, the probability of these events does not depend on $x$. We denote it by $\pi(\alpha)$:

$$\pi(\alpha) = P(G(0, \alpha)).$$

To deal with large radii, we introduce two other families of events as follows. For all $\alpha > 0$, we define an event $H(\alpha)$ by

$$H(\alpha) = \{\exists c \in \chi \setminus B(0, 10\alpha) : B(c, r(c)) \cap B(0, 9\alpha) \neq \varnothing\}$$

and an event $\widetilde{H}(\alpha)$ by

$$\widetilde{H}(\alpha) = \{\exists c \in \chi \cap B(0, 100\alpha) : r(c) \geq \alpha\}.$$

Finally, we define a random variable $M$ as follows:

$$M = \sup_{x \in S} \|x\|. \tag{4}$$

(We let $M = 0$ if $S$ is empty.)

Our aim in this subsection is to prove the following proposition.

PROPOSITION 3.1. *There exists a constant $C$ that depends only on the dimension $d$, such that the following assertions hold for all $\alpha > 0$:*

$$\pi(10\alpha) \leq C\pi(\alpha)^2 + \lambda C \int_\alpha^{+\infty} r^d \nu(dr), \tag{5}$$

$$P(M \geq 9\alpha) \leq \pi(\alpha) + \lambda C \int_\alpha^{+\infty} r^d \nu(dr) \tag{6}$$

*and*

$$\pi(\alpha) \leq C\lambda \alpha^d. \tag{7}$$

In the following lemma, we provide a link between the percolation event and the families of events $G(0,\cdot)$ and $H(\cdot)$. This gives part of the proof of (6).

LEMMA 3.2. *For all $\alpha > 0$, the following inclusion holds:*
$$\{M \geq 9\alpha\} \subset G(0,\alpha) \cup H(\alpha).$$

PROOF. Let $\alpha > 0$. If $G(0,\alpha)$ does not occur, then one cannot go from 0 to the complement of $B(0, 8\alpha)$ using balls $B(c, r(c))$, $c \in \chi \cap B(0, 10\alpha)$. If, moreover, $H(\alpha)$ does not occur, then the balls $B(c, r(c))$, $c \in \chi \setminus B(0, 10\alpha)$, will not help to connect 0 to the complement of $B(0, 8\alpha)$. Therefore, one has $M \leq 8\alpha$. □

In the following lemma, we provide a means to control the probabilities $\pi(\alpha)$. This supplies part of the proof of (5).

LEMMA 3.3. *There exists a constant $C_1$ that depends only on the dimension $d$ such that, for all $\alpha > 0$, the following holds:*
$$\pi(10\alpha) \leq C_1 \pi(\alpha)^2 + P(\widetilde{H}(\alpha)).$$

PROOF. For all $r \geq 0$, we denote by $S_r$ the Euclidean sphere centered at the origin with radius $r$:
$$S_r = \{x \in \mathbb{R}^d : \|x\| = r\}.$$

Fix $K$ and $L$, two subsets of $\mathbb{R}^d$ such that the following properties hold:

1. the sets $K$ and $L$ are finite;
2. $K \subset S_{10}$ and $L \subset S_{80}$;
3. $S_{10} \subset K + B(0,1)$ and $S_{80} \subset L + B(0,1)$.

We define $C_1$ as the product of the cardinalities of the sets $K$ and $L$.

Let $\alpha > 0$. In this step, we prove the following inclusion:

$$(8) \qquad G(0, 10\alpha) \setminus \widetilde{H}(\alpha) \subset \left( \bigcup_{k \in K} G(\alpha k, \alpha) \right) \cap \left( \bigcup_{l \in L} G(\alpha l, \alpha) \right).$$

We assume that the event $G(0, 10\alpha)$ occurs, but that the event $\widetilde{H}(\alpha)$ does not occur. As $G(0, 10\alpha)$ occurs, one can go from $S_{10\alpha}$ to $S_{80\alpha}$ using only balls $B(c, r(c))$ centered at points of $\chi \cap B(0, 100\alpha)$. We refer to Figure 1 where the dotted circles stand for some of the previous balls $B(c, r(c))$. One of these balls touches $S_{10\alpha}$. This ball then touches $B(\alpha k, \alpha)$ for some $k \in K$. We then see that one can go from $B(\alpha k, \alpha)$ to the complement of $B(\alpha k, 8\alpha)$ using only balls $B(c, r(c))$ centered at points of $\chi \cap B(0, 100\alpha)$. But, as $\widetilde{H}(\alpha)$ does not occur, the radius of each such ball $B(c, r(c))$ is less than $\alpha$. Therefore, one can go from $B(\alpha k, \alpha)$ to the complement of $B(\alpha k, 8\alpha)$ using only balls $B(c, r(c))$ centered at points of $\chi \cap B(\alpha k, 10\alpha)$. In other words, $G(\alpha k, \alpha)$ occurs. We have proven that the event $\bigcup_{k \in K} G(\alpha k, \alpha)$ occurs. We can prove in a similar way that the event $\bigcup_{l \in L} G(\alpha l, \alpha)$ occurs. Therefore the inclusion (8) is proved.

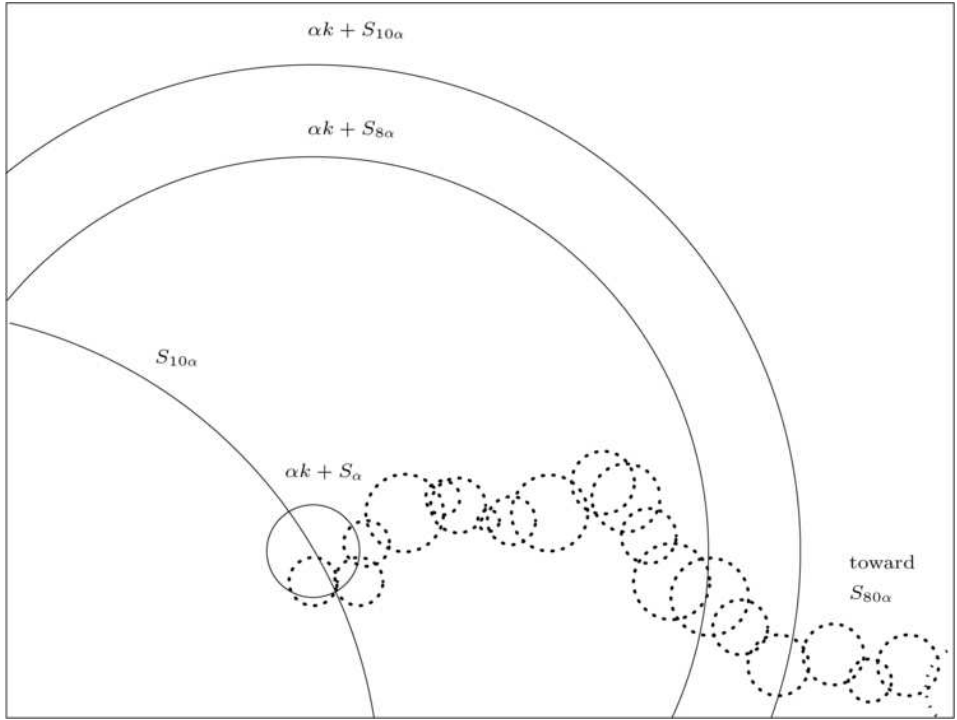

Fig. 1. *Proof of* (8).

The right-hand side of (8) is an intersection of two events. The first of them only depends on what occurs in $B(0, 20\alpha)$. The other one only depends on what occurs in $B(0, 70\alpha)^c$. These two events are therefore independent. The result then follows from relation (8). $\square$

In the following two lemmas, we provide a means to bound the probabilities of the events $P(H(\alpha))$ and $P(\widetilde{H}(\alpha))$. This will enable us to conclude the proof of (5) and (6).

LEMMA 3.4. *There exists a constant $C_2$ that depends only on the dimension $d$, such that for all $\alpha > 0$, the following inequality holds:*

$$P(H(\alpha)) \leq \lambda C_2 \int_\alpha^{+\infty} r^d \nu(dr).$$

PROOF. Let $\alpha > 0$. We have

$$H(\alpha) = \{\xi \cap V(\alpha) \neq \varnothing\},$$

where

$$V(\alpha) = \{(c,r) \in \mathbb{R}^d \times ]0, +\infty[ : c \notin B(0, 10\alpha) \text{ and } B(c,r) \cap B(0, 9\alpha) \neq \varnothing\}.$$

We therefore have

$$\begin{aligned} P(H(\alpha)) &= P(\xi \cap V(\alpha) \neq \varnothing) \\ &\leq E(\text{card}(\xi \cap V(\alpha))) \\ &= \int_{\mathbb{R}^d} \lambda \, dc \int_0^{+\infty} \nu(dr) 1_{V(\alpha)}(c,r). \end{aligned}$$

As

$$V(\alpha) = \{(c,r) \in \mathbb{R}^d \times ]0, +\infty[ : c \notin B(0, 10\alpha) \text{ and } c \in B(0, 9\alpha + r)\},$$

we obtain

$$\begin{aligned} P(H(\alpha)) &\leq \lambda \int_0^{+\infty} \nu(dr) |B(0, 9\alpha + r) \setminus B(0, 10\alpha)| \\ &= \lambda \int_\alpha^{+\infty} |B(0, r + 9\alpha)| \nu(dr) \\ &\leq \lambda \int_\alpha^{+\infty} |B(0, 10r)| \nu(dr). \end{aligned}$$

The inequality stated in the lemma is therefore fulfilled with $C_2 = |B(0,10)|$. □

LEMMA 3.5. *There exists a constant $C_3$ that depends only on the dimension $d$, such that for all $\alpha > 0$, the following inequality holds:*

$$P(\widetilde{H}(\alpha)) \leq \lambda C_3 \int_\alpha^{+\infty} r^d \nu(dr).$$

PROOF. Let $\alpha > 0$. We have

$$\begin{aligned} P(\widetilde{H}(\alpha)) &\leq E(\text{card}(\{c \in \chi \cap B(0, 100\alpha) : r(c) \geq \alpha\})) \\ &= \lambda |B(0, 100\alpha)| \nu([\alpha, +\infty[) \end{aligned}$$

$$= \lambda |B(0,100)| \int_\alpha^{+\infty} \alpha^d \nu(dr)$$

$$\leq \lambda |B(0,100)| \int_\alpha^{+\infty} r^d \nu(dr).$$

The inequality stated in the lemma is therefore fulfilled with $C_3 = |B(0,100)|$. □

The following lemma will enable us to make sure that $\pi$ is small enough on a sufficiently large set. This will give (7).

LEMMA 3.6. *There exists a constant $C_4$ that depends only on the dimension $d$, such that for all $\alpha > 0$, the following inequality holds:*

$$\pi(\alpha) \leq \lambda C_4 \alpha^d.$$

PROOF. Let $\alpha > 0$. Note that, if $B(0,10\alpha) \cap \chi$ is empty, then $\Sigma(B(0,10\alpha))$ is empty and therefore the event $G(0,\alpha)$ cannot occur. As a consequence,

$$P(G(0,\alpha)) \leq P(B(0,10\alpha) \cap \chi \neq \varnothing)$$
$$\leq E(\operatorname{card}(B(0,10\alpha) \cap \chi))$$
$$= \lambda |B(0,10\alpha)|.$$

The inequality stated in the lemma is therefore satisfied with $C_4 = |B(0,10)|$. □

PROOF OF PROPOSITION 3.1.  This is a consequence of Lemmas 3.2, 3.3, 3.4, 3.5 and 3.6.  □

3.2. *Proof of the existence of subcritical behavior.* We need the following lemma.

LEMMA 3.7. *Let $f$ and $g$ be two measurable, bounded and nonnegative functions from $[1,+\infty]$ to $\mathbb{R}_+$. We assume that $f$ is bounded by $1/2$ on $[1,10]$ and that $g$ is bounded by $1/4$ on $[1,+\infty]$. We also assume that, for all real $\alpha \geq 10$, the following inequality holds:*

(9) $$f(\alpha) \leq f(\alpha/10)^2 + g(\alpha).$$

*Under those conditions, if $g(\alpha)$ converges to $0$ as $\alpha$ tends to infinity, then $f(\alpha)$ converges to $0$ as $\alpha$ tends to infinity. If, moreover, a real number $s$ is such that the integral $\int_1^{+\infty} \alpha^s g(\alpha)\,d\alpha$ is finite, then the integral $\int_1^{+\infty} \alpha^s f(\alpha)\,d\alpha$ is also finite.*

PROOF. We assume that $g$ converges to 0. As $f$ is bounded by $1/2$ on $[1,10]$ and $g$ is bounded by $1/4$ on $[1,+\infty]$, we obtain, by (9), that $f$ is bounded by $1/2$ on $[1,+\infty]$. Therefore, for all real $\alpha \geq 10$, we have

$$f(\alpha) \leq f(\alpha/10)/2 + g(\alpha).$$

As a consequence, for all $\alpha \in [1,10]$ and all integers $n \geq 1$, the following inequality holds:

$$\begin{aligned}
f(10^n\alpha) &\leq f(\alpha)/2^n + g(10\alpha)/2^{n-1} + \cdots + g(10^n\alpha) \\
&\leq 1/2^{n+1} + g(10\alpha)/2^{n-1} + \cdots + g(10^n\alpha).
\end{aligned} \tag{10}$$

For all integers $n \geq 1$, we let

$$F_n = \sup_{\alpha \in [1,10]} f(10^n\alpha) \quad \text{and} \quad G_n = \sup_{\alpha \in [1,10]} g(10^n\alpha).$$

By (10), we obtain

$$F_n \leq 1/2^{n+1} + G_1/2^{n-1} + \cdots + G_n. \tag{11}$$

As $g$ is bounded and converges to 0, the sequence $(G_n)_n$ converges to 0. By (11), we then obtain the convergence of $F_n$ to 0. Therefore, as $f$ is nonnegative, $f$ converges to 0.

Let $s$ be a real number. We furthermore assume that the integral $\int_1^{+\infty} \alpha^s g(\alpha)\,d\alpha$ is finite. By the first step, we know that $f$ converges to 0. Therefore, there exists a real $A \geq 10$ that we fix, such that $f(\alpha)$ is bounded by $10^{-s-1}/2$ on $[A/10, +\infty[$. For all real $r \geq A$, we obtain, by (9),

$$\begin{aligned}
\int_A^r f(\alpha)\alpha^s\,d\alpha &\leq \int_A^r f(\alpha/10)^2 \alpha^s\,d\alpha + \int_A^r g(\alpha)\alpha^s\,d\alpha \\
&\leq 10^{s+1} \int_{A/10}^{r/10} f(\alpha)^2 \alpha^s\,d\alpha + \int_A^{+\infty} g(\alpha)\alpha^s\,d\alpha \\
&\leq 1/2 \int_{A/10}^{r/10} f(\alpha)\alpha^s\,d\alpha + \int_A^{+\infty} g(\alpha)\alpha^s\,d\alpha \\
&\leq 1/2 \int_A^r f(\alpha)\alpha^s\,d\alpha + 1/2 \int_{A/10}^A f(\alpha)\alpha^s\,d\alpha + \int_A^{+\infty} g(\alpha)\alpha^s\,d\alpha.
\end{aligned}$$

We therefore obtain

$$\int_A^r f(\alpha)\alpha^s\,d\alpha \leq \int_{A/10}^A f(\alpha)\alpha^s\,d\alpha + 2\int_A^{+\infty} g(\alpha)\alpha^s\,d\alpha.$$

As a consequence, the integral $\int_A^{+\infty} f(\alpha)\alpha^s\,d\alpha$ is finite. $\square$

The following result gives one direction in the equivalences stated in the theorems.

LEMMA 3.8. *There exists a positive constant $\widetilde{C}$ that depends only on the dimension $d$, such that the following assertions hold for all probability measure $\nu$ on $]0, +\infty[$:*

1. *if $E(R^d)$ is finite, then $P_{\lambda,\nu}(Percolation) = 0$ for all positive real $\lambda < \widetilde{C}(E(R^d))^{-1}$;*
2. *for all $s > 0$, if $E(R^{d+s})$ is finite, then $E_{\lambda,\nu}(M^s)$ is finite for all positive real $\lambda < \widetilde{C}(E(R^d))^{-1}$.*

PROOF. Let $C$ be the constant given by Proposition 3.1. We define a constant $\widetilde{C}$ by

$$\widetilde{C} = (4C^2)^{-1}.$$

Let $\nu$ be a probability measure on $]0, +\infty[$. Throughout the proof, we assume that $E(R^d)$ is finite. We set

$$\lambda_0 = \widetilde{C}(E(R^d))^{-1}.$$

Let $\lambda > 0$ be such that $\lambda < \lambda_0$.

Let us define a positive real number $A$ by

$$A = (E(R^d))^{1/d}/10.$$

Let $f: [1, +\infty[$ be the function defined by

$$f(\alpha) = C\pi(A\alpha)$$

and $g: [1, +\infty[$ be the function defined by

$$g(\alpha) = \lambda C^2 \int_{A\alpha/10}^{\infty} r^d \nu(dr).$$

As $\lambda < \lambda_0$, we obtain, by (7), that $f$ is bounded by $1/2$ on $[1, 10]$. As $\lambda < \lambda_0$, we obtain that $g$ is bounded by $1/4$ on $[1, +\infty[$. As $E(R^d)$ is finite, the function $g$ converges to 0. By (5), we obtain that (9) holds for all $\alpha \geq 10$. By Lemma 3.7 we therefore obtain that $f$, and then $\pi$, converges to 0. By (6), we then obtain that $M$ is almost surely finite. Therefore, almost surely, percolation does not occur.

Let $s > 0$. In this step, we assume, furthermore, that $E(R^{d+s})$ is finite. The integral $\int_1^{+\infty} \alpha^{s-1} g(\alpha) \, d\alpha$ is therefore finite. By Lemma 3.7, we obtain that the integral $\int_1^{+\infty} \alpha^{s-1} f(\alpha) \, d\alpha$ is also finite. By (6), we then obtain that the integral $\int_1^{+\infty} \alpha^{s-1} P_{\lambda,\nu}(M \geq 9A\alpha) \, d\alpha$ is also finite. As a consequence, the moment $E_{\lambda,\nu}(M^s)$ is finite. □

3.3. *Proof of nonexistence of subcritical behavior.* In the following lemma, we supply the other direction of the equivalences stated in the theorems. Recall that $\Sigma$ is defined by (1) and that $M$ is defined by (4). Let us recall that the first part of this lemma is a special case of a result by Meester and Roy (see [5], Proposition 7.3).

LEMMA 3.9. *Let $\nu$ be a probability measure on $]0,+\infty[$. If $E(R^d)$ is infinite then, for all $\lambda > 0$, we have $P_{\lambda,\nu}$-almost surely $\Sigma = \mathbb{R}^d$. If $s > 0$ is such that $E(R^{d+s})$ is infinite, then for all $\lambda > 0$, $E_{\lambda,\nu}(M^s)$ is infinite.*

PROOF. Let $\nu$ be a probability measure on $]0,+\infty[$ and $\lambda > 0$. We first prove that for all $r > 0$, the following inequality holds:

$$
\begin{aligned}
P_{\lambda,\nu}(\exists c \in \chi : B(0,r) \subset B(c,r(c))) \\
\geq 1 - \exp\biggl(-\lambda 2^{-d}|B(0,1)| \int_{[2r,+\infty[} \alpha^d \nu(d\alpha)\biggr).
\end{aligned}
\tag{12}
$$

Let $r > 0$. We have

$$
\begin{aligned}
P_{\lambda,\nu}(\exists c \in \chi : B(0,r) \subset B(c,r(c))) &= 1 - \exp\biggl(-\lambda \int_{\mathbb{R}^d} P(R \geq \|x\| + r)\, dx\biggr) \\
&= 1 - \exp(-\lambda E(|B(0,R-r)|1_{R \geq r})) \\
&\geq 1 - \exp(-\lambda E(|B(0,R-r)|1_{R \geq 2r})) \\
&\geq 1 - \exp(-\lambda E(|B(0,R/2)|1_{R \geq 2r})).
\end{aligned}
$$

Relation (12) is thus proved.

If $E(R^d)$ is infinite, then, by (12), we obtain, for all $r > 0$,

$$P_{\lambda,\nu}(\exists c \in \chi : B(0,r) \subset B(c,r(c))) = 1.$$

Therefore, almost surely, we have $\Sigma = \mathbb{R}^d$.

Let $s > 0$. We now assume that $E(R^{d+s})$ is infinite. If $E(R^d)$ is infinite, then the desired result is a consequence of what we have proven in the previous step. We assume, henceforth, that $E(R^d)$ is finite. Let $C$ be defined by

$$C = \lambda 2^{-d}|B(0,1)| \int_{[0,+\infty[} \alpha^d \nu(d\alpha).$$

This constant is finite. By (12), we obtain, for all $r > 0$, the following inequality:

$$
\begin{aligned}
P_{\lambda,\nu}(\exists c \in \chi : B(0,r) \subset B(c,r(c))) \\
\geq C^{-1}(1 - \exp(-C))\lambda 2^{-d}|B(0,1)| \int_{[2r,+\infty[} \alpha^d \nu(d\alpha).
\end{aligned}
$$

As $E(R^{d+s})$ is infinite, the integral

$$\int_0^{+\infty} \left( r^{s-1} \int_{2r}^{+\infty} \alpha^d \nu(d\alpha) \right) dr$$

is infinite. Therefore, the integral

$$\int_0^{+\infty} r^{s-1} P_{\lambda,\nu}(\exists c \in \chi : B(0,r) \subset B(c, r(c))) \, dr$$

is infinite. As a consequence, the integral $\int_0^{+\infty} r^{s-1} P_{\lambda,\nu}(M \geq r) \, dr$ is infinite. The moment $E_{\lambda,\nu}(M^s)$ is then infinite. $\square$

REMARK. Using Lemma 3.8 and the proof of Lemma 3.9, we could also give necessary and sufficient conditions for the integrability of the volume of $S$, or for the integrability of the radius of the largest ball centered at the origin and contained in $S$.

**Acknowledgment.** I would like to thank the referee for careful reading and helpful remarks.

## REFERENCES

[1] BŁASZCZYSZYN, B., RAU, C. and SCHMIDT, V. (1999). Bounds for clump size characteristics in the Boolean model. *Adv. in Appl. Probab.* **31** 910–928. MR1747448
[2] HALL, P. (1985). On continuum percolation. *Ann. Probab.* **13** 1250–1266. MR0806222
[3] HALL, P. (1988). *Introduction to the Theory of Coverage Processes.* Wiley, New York. MR0973404
[4] KALLENBERG, O. (1986). *Random Measures*, 4th ed. Akademie-Verlag, Berlin. MR0854102
[5] MEESTER, R. and ROY, R. (1996). *Continuum Percolation.* Cambridge Univ. Press. MR1409145
[6] MENSHIKOV, M. V., POPOV, S. YU. and VACHKOVSKAIA, M. (2003). On a multiscale continuous percolation model with unbounded defects. *Bull. Braz. Math. Soc. (N.S.)* **34** 417–435. MR2045167
[7] MØLLER, J. (1994). *Lectures on Random Voronoĭ Tessellations.* Springer, New York. MR1295245
[8] NEVEU, J. (1977). Processus ponctuels. *École d'Été de Probabilités de Saint-Flour VI—1976. Lecture Notes in Math.* **598** 249–445. Springer, Berlin. MR0474493

UNIVERSITÉ D'ORLÉANS
MAPMO
B.P. 6759
45067 ORLÉANS CEDEX 2
FRANCE
E-MAIL: jbgouere@univ-orleans.fr